\documentclass[a4paper,12pt]{article}
\usepackage[dvips]{graphicx}
\usepackage{color}
\usepackage{graphicx}
\usepackage{amsmath,amssymb,amsthm}
\usepackage{wrapfig}

\newtheorem{thm}{Theorem}
\newtheorem{cor}{Corollary}[section]
\newtheorem{prop}{Proposition}[section]
\newcounter{def}

\newcommand{\cm}{{\cal M}}
\newcommand{\dist}{\mathop{\rm dist}\nolimits}
\newcommand{\cg}{\mathop{\rm CG}\nolimits}

\begin{document}

\title{Enumeration of irreducible  contact graphs on the sphere}
\author{Oleg R. Musin\thanks{This research is supported by the Russian government project 11.G34.31.0053,
  RFBR grant 11-01-00735, and NSF grant DMS 1101688.} \, and Alexey S. Tarasov\thanks{	This research is supported by the Russian
    government project 11.G34.31.0053 and 
    RFBR grant 11-01-00735.}}
\date{}
\maketitle
\begin{abstract}
In this article, using the computer, are enumerated all locally-rigid packings by $N$ congruent circles (spherical caps) on the unit sphere ${\Bbb S}^2 $ with $N < 12.$ This is equivalent to the enumeration of irreducible spherical contact graphs.
\end{abstract}

\section{Introduction}

Packings where all spheres are constrained by their neighbors to stay in one location are called {\it rigid} or {\it locally-rigid.} So every sphere of this packing jammed by neighbors and it can not be shifted to the side in order to increase the minimum distance between the center of the sphere and the other centers of the balls.  

Consider the $N$ non-overlapping spheres of the same radius $r$ in $ {\Bbb R}^ 3$ , which are arranged so that they all touch  to one (central) sphere of unit radius. We denote by $P: = \{A_{1}, \ldots, A_{N}\}$ set of points where external spheres touch the central sphere. Join the points $A_ {i}$ and $A_ {j}$ by edge (minimum arc of a great circle  if the relevant external spheres touch. The resulting graph is called {\it contact} and denote $\cg(P)$. If this packing on ${\Bbb S}^2$  is a locally rigid , then we say that the graph $\cg(P)$ {\it irreducible}.
Thus, the problem  of studying locally rigid packings reduces to the study of irreducible contact graphs.

There are several connections between this geometric problem with other  sphere packings problems. The main application outside mathematics is ``jammed'' (locally rigid) hard-particles packings in Materials Science   
(see, for instance, \cite{Appl04, Appl10}). Note that for most potentials in Physics minimum energy configurations of particles are also locally rigid. 

In Mathematics W. Habicht, K. Sch\"utte,  B.L. van der Waerden, and L. Danzer  applied irreducible contact graphs  for the kissing number  and   Tammes problems \cite{HabvdW, SvdW1, vdW, SvdW2, Dan}. Let us consider briefly these two classical geometric problem. 

\medskip

The {\it kissing number} $k_3$ is the highest number of equal
non-overlapping spheres in ${\Bbb R}^3$ that touch another
sphere of the same size. In other words, the kissing number
problem asks how many white billiard balls can
{\em kiss} (touch) a black ball.

The most symmetrical configuration, 12 balls around
another, is achieved if the 12 balls are placed at positions
corresponding to the vertices of a regular icosahedron concentric
with the central ball. However, these 12 outer balls do not kiss
each other and may all be moved freely. This space between the balls introduces a question: If you moved
all of them to one side, would a 13th ball fit?

This problem was the subject of the famous discussion between Isaac
Newton and David Gregory in 1694. Most reports say that Newton
believed the answer was 12 balls, while Gregory thought that 13 might be possible. That  is why it often called the {\it thirteen spheres problem}

The problem was finally solved by Sch\"utte and van der Waerden
in 1953 \cite{SvdW2}. A subsequent two-page sketch of an elegant
proof was given by Leech \cite{Lee} in 1956.
%Most people agree that Leech's proof is correct, but there are gaps in his exposition,
%many of them involving sophisticated spherical trigonometry.
Leech's proof was presented in the first edition  of the well-known book
by Aigner and Ziegler \cite{AZ}; the authors removed this chapter from the second edition because a complete proof would have to include so much spherical trigonometry.

The thirteen spheres problem continues to be of interest, and
new proofs have been published in the last several years by Hsiang
\cite{Hs}, Maehara \cite{Ma, Manew} (this proof is based on Leech's proof),  B\"or\"oczky \cite{Bor},
Anstreicher \cite{Ans}, and Musin \cite{Mus13}. 

\medskip

If $N$ unit spheres kiss the unit sphere in ${\Bbb R}^3$, then
the set of kissing points
is an arrangement on the central sphere such that the (Euclidean)
distance between any two points is at least 1. This allows us to state the kissing number problem in another way: How many points can
be placed on the surface of ${\Bbb S}^2$ so that the angular
separation between any two points be at least $60^{\circ}$?

This leads to an important generalization: to find a set $X$ of $N$ point in ${\Bbb S}^{2}$  such that the minimum angular distance of distinct points in $X$ is large as possible. In other words,
{ how are $N$ congruent, non-overlapping circles distributed on the sphere
when the common radius of the circles has to
be as large as possible?}

%his question, also known as the problem of the ``inimical dictators'': {\it Where should $N$ dictators build their palaces on a planet so as to be as far away from each other as possible?} 

The problem
was first asked by the Dutch botanist Tammes \cite{Tam} (see \cite[Section 1.6: Problem 6]{BMP}), who was led to this problem by examining the distribution of
openings  on the pollen grains of different flowers.

The Tammes problem is presently solved only for several values
of $N$: for $N=3,4,6,12$ by L. Fejes T\'oth \cite{FeT0}; for
$N=5,7,8,9$ by Sch\"utte
and van der Waerden \cite{SvdW1}; for $N=10,11$ by Danzer \cite{Dan} (for  $N=11$ see also B\"or\"oczky \cite{Bor11}) and for $N=24$ by Robinson \cite{Rob}. We recently solved Tammes' problem for the case $N=13$  \cite{MT}.

Note that the kissing number problem currently solved only for dimensions  $n=3,4,8$  and $24$ (see \cite{BDM, Mus1, Mus3, Mus4}). Proofs in these papers are based on the Delsarte method and its generalizations (see, for example,  \cite{BM, CS, Mus5, Mus6}).

\section{The irreducible contact graphs}

\subsection{Basic definitions}
Let $X$ be a finite subset of ${\Bbb S}^{2}$. Denote
$$\psi(X):=\min\limits_{x,y\in X}{\{\dist(x,y)\}}, \mbox{ where } x\ne y.$$
%Then $X$ is a spherical $\psi(X)$-code.

Denote by $d_N$ the largest angular separation $\psi(X)$ with $|X|=N$ that can be attained
 in ${\Bbb S}^{2}$, i.e.
$$
d_N:=\max\limits_{X\subset{\Bbb S}^2}{\{\psi(X)\}}, \, \mbox{ where } \;  |X|=N.
$$

\noindent{\bf Contact graphs.} Let $X$ be a finite set in ${\Bbb
S}^2$. The {\it contact graph} $\cg(X)$  is the graph with vertices
in $X$ and edges $(x,y), \, x,y\in X$ such that $\dist(x,y)=\psi(X)$.

\medskip

\noindent{\bf Shift of a single vertex.}  
We say that a vertex $x\in X$  {\it can be shifted},  if in any open neighborhood of $x$  there is a point  $x'\in {\Bbb S}^2$ such that $$\dist(x',X\setminus\{x\})>\dist(x,X\setminus\{x\}).$$

\medskip

\noindent{\bf Irreducible graphs.}
We say that the graph $\cg(X)$ is {\it irreducible} if there are  no shifts of vertices. 
This terminology was used by Sch\"utte - van der Waerden \cite{SvdW1,SvdW2}, Fejes T\'oth \cite{FeT}, and Danzer \cite{Dan}.

Let us denote by $J_N$ the family of all sets $X$ in ${\Bbb S}^2$ such that $|X|=N$ and its contact graph $\cg(X)$ is irreducible.

\medskip

\noindent{\bf  D-irreducible graphs.}
L. Danzer \cite[Sec. 1]{Dan} defined the following move of a vertex. Let  
$x,y,z$ be vertices of $\cg(X)$ with $\dist(x,y)=\dist(x,z)=\psi(X)$. Denote by  $x^0$ the mirror image of $x$ with respect 
the great circle passes through  $y, z$ (see Fig~\ref{fig3}). We call this move as 
 {\it D-flip}, if $\dist(x^0,X\setminus\{x,y,z\}) > \psi(X)$.

An irreducible contact graph $\cg(X)$ is called D-irreducible if it does not admit any D-flip.    

\medskip

\begin{figure}[h]
\begin{center}
\includegraphics[clip,scale=1]{pics/13-6.mps}
\end{center}
\caption{D-flip}
\label{fig3}
\end{figure}

\medskip

\noindent{\bf Maximal graphs.} Suppose that $X\subset{\Bbb S}^{2}$  with $|X|=N$ and $\psi(X)=d_N$.  Then we call this contact graph $\cg(X)$ - {\it maximal.}

\subsection{Properties of irreducible contact graphs}

In this subsection we consider $X\subset {\Bbb S}^2$ such that the graph $\cg(X)$ is irreducible, i.e. $X\in J_N$.     
%of all  In other words, $X\in J_N$, where $N=|X|$. The lengths of all edges of $\cg(X)$ equal $\psi(X)$. 
The following properties of  $J_N$  were found in \cite{SvdW1}, \cite{Dan}, and \cite{BS,BS14} (see also  \cite[Chap. VI]{FeT}). 

 Let $a,b,x,y\in X$ with $\dist(a,b)=\dist(x,y)=\psi(X)$.
Then the shortest arcs ${ab}$ and ${xy}$ don't intersect. Otherwise,
the length of at least one of the arcs $ax, ay, bx, by$ has to be
less than $\psi(X)$. This yields the planarity of $\cg(X)$.
\begin{prop} If $X$ is a finite subset of ${\Bbb S}^2$, then  $\cg(X)$ is a planar graph.
\end{prop}

\begin{prop} If $X\in J_N$, then all faces of $\cg(X)$ are convex polygons in  ${\Bbb S}^2$. 
\end{prop}
(Indeed, otherwise, a ``concave'' vertex of a face $P$ can be shifted to the interior of $P$.) 

\medskip

 Let $X$ be a subset of ${\Bbb S}^2$ with $|X|=N$.  We say that $X$ is {\it maximal} if $\psi(X)=d_N$.

\begin{prop} If $X$ is maximal, then for  $N>5$ the graph $\cg(X)$ is irreducible.
\end{prop}

\begin{prop} If $X\in J_N$, then degrees of its vertices can take only the values $0$ (isolated vertices), $3$, $4$, or  $5$.
\end{prop}

\begin{prop} If $X\in J_N$, then  faces of $\cg(X)$ are polygons with at most $\lfloor2\pi/\psi(X)\rfloor$ vertices.
\end{prop}

The following simple proposition has been proved by	B\"or\"oczky and Szab\'o in  {\cite[Lemma 8 and Lemma 9(iii)]{BS}}. Actually,  they considered the case $N=13$. However, the proof works for all $N$. 

\begin{prop}  Let $X\in J_N$. If  $\cg(X)$ contains an isolated vertex, then it lies in the interior of a polygon of $\cg(X)$ with six or more vertices. Moreover, if it is a hexagon, then it cannot contain two isolated vertices. 
\end{prop}

Combining these propositions, we obtain the following combinatorial properties of irreducible contact graphs.
\begin{cor} If $X\in J_N$, then  $G:=\cg(X)$ satisfies the following properties 
\begin{enumerate}
\item $G$ is a planar graph;
\item Any vertex of $G$ is of degree $0,3,4,$ or $5$;
\item If  $G$ contains an isolated vertex $v$, then $v$ lies in a  face with $m\ge 6$ vertices. Moreover, a hexagonal face of  $G$ cannot contain two or more isolated vertices.
\label{cor1}
\end{enumerate}

\end{cor}

\section{Danzer's work on irreducible contact graphs} 

L. Danzer \cite{Dan} solved the Tammes problem for $N=10$ and $N=11$. His proof is based on the concept of irreducible graphs. (Actually, this paper is a translation of the Habilitationsschrift of Ludwig Danzer ``Endliche Punktmengen auf der
2-sph\"are mit m\"oglichst gro{\ss}em Minimalabstand''. Universit\"at G\"othingen, 1963) In particular,  he added to shifting a single vertex a new idea - a shift which we call here {D--flip,} i. e. {\it Danzer's flip.}

\medskip

In \cite{Dan} Danzer gives the list of all D-irreducible graphs for $6\le N \le 10$. Since the contact graph of a maximal set is irreducible (and D-irreducible \cite{Dan}), this list implies a solution of the Tammes problem for  $6\le N \le 10$. (For the case $N=11$ Danzer considered only maximal sets.) 

Here we give the Danzer list of D-irreducible graphs.

\subsection{$N=6$}

\begin{tabular}{cccc}
maximal  &
\includegraphics[clip,scale=0.2]{pics/irr-4.mps}  
&  $ G_I = \cm_6(t)$ &
\includegraphics[clip,scale=0.3]{pics/irr-6.mps}  
\end{tabular}

\subsection{$N=7$}

\begin{tabular}{cccc}
maximal $\cm_7 = G_{III}$ & 
\includegraphics[clip,scale=0.5]{pics/seven2.mps} 
&  $\cm_7(t)=G_{II}$&
\includegraphics[clip,scale=0.5]{pics/seven1.mps} 
\end{tabular}

\subsection{$N=8$}

\begin{tabular}{cccc}
\includegraphics[clip,scale=0.5]{pics/eight4.mps} 
& maximal $\cm_8 = G_V$ & & \\
\includegraphics[clip,scale=0.5]{pics/eight2.mps} 
& $\cm_8(t)=G_{IV}$ &
\includegraphics[clip,scale=0.5]{pics/eight3.mps} 
& $\cm_8(u,v)=G_{IX}$ 
\end{tabular}

\subsection{$N=9$}

\begin{tabular}{cccc}
\includegraphics[clip,scale=0.5]{pics/nine7.mps} 
& maximal & & \\К
\includegraphics[clip,scale=0.5]{pics/nine2.mps} 
& $\cm_9(t)$ & 
\includegraphics[clip,scale=0.5]{pics/nine9.mps} 
& $\cm_9^*$ 

\end{tabular}

\subsection{$N=10$}

\begin{tabular}{cccc}
\includegraphics[clip,scale=0.5]{pics/ten17.mps} 
& maximal $\cm_{10}=G_{XI}$ & & \\
\includegraphics[clip,scale=0.5]{pics/ten21.mps} 
& $\cm_{10,2}(t) = G_{XII}$ &
\includegraphics[clip,scale=0.5]{pics/ten27.mps} 
& $\cm_{10,3}(t) = G_{XIII}$ \\
\includegraphics[clip,scale=0.5]{pics/ten29.mps} 
& $\cm^{*} = G_{XVIII}$ & 
\includegraphics[clip,scale=0.5]{pics/ten28.mps} 
& $\cm^{**} = G_{XIX}$ \\
\includegraphics[clip,scale=0.5]{pics/ten10.mps} 
& $\cm^3_{10}(t) = G_{XIV}$ & 
\includegraphics[clip,scale=0.5]{pics/ten7.mps} 
& $\cm^{1,2}_{10}(t)=G_{XV}$ \\
\includegraphics[clip,scale=0.5]{pics/ten25.mps} 
& $\cm^{1,3}_{10}(t)=G_{XVI}$ &
\includegraphics[clip,scale=0.5]{pics/iv_ten1.mps} 
& $\tilde \cm^{1,2,3}_{10}(t) \cup {\cm} = G_{VIII}$ \\
\includegraphics[clip,scale=0.5]{pics/ten23.mps} 
& $\cm^{3}_{10}(u,v)=G_X$ & & 
\end{tabular}

\section{Enumeration of irreducible contact graphs.}

\subsection{Geometric embedding of irreducible contact graphs.}

Let $X\subset {\Bbb S}^2$ be a finite point set such that  its contact graph $\cg(X)$ is irreducible. In Corollary 2.1 we collected together combinatorial properties of $\cg(X)$.  There are several geometric properties. 

Recall that all faces of  $\cg(X)$ are convex (Proposition  2.2). Since all edges of $\cg(X)$ have the same lengths $\psi(X)$, all its faces are spherical equilateral convex polygon with number of vertices at most $\lfloor2\pi/\psi(X)\rfloor$.

Рассмотрим теперь планарный граф  Consider now a planar graph $G$ with given faces  $\{F_k\}$ that satisfy Corollary 2.1. We are going consider embeddings of this graph into  ${\Bbb S}^2$ as an irreducible contact graph  $\cg(X)$ for some $X\subset {\Bbb S}^2$.

Any embedding of $G$ in ${\Bbb S}^2$ is uniquely defined by the following list of parameters (variables):\\
(i) The edge length  $d$;\\
(ii) The set of all angles $u_{ki}$, $i=1,\ldots,m_k$ of faces $F_k$. Here is $m_k$ denotes the number of vertices of $F_k$.)

In our paper where we give a solution of Tammes' problem for  $N=13$ \cite{MT} was considered main relations between these parameters  (\cite[Propositions 3.6--3.11]{MT}). Let us give here these results.  (We added also a general statement for $m>4$.)

 %Since equations and correspondent linear inequalities for faces with five and especially with six vertices are rather complicated for some graphs our program worked a lot of time. Now we decided triangulate all faces with $m>3$ vertices.   

% В частности, для всех $u_{ki}$ верны следующие утверждения. 

\begin{prop}
\begin{enumerate}
\item $u_{ki} <\pi$ for all $i$ and $k$.
\item $u_{ki} \geqslant \alpha(d)$ for all $i$ and $k$, where 
%$\alpha:=\arccos(\cos{d}/(1+\cos{d}))$
$$ \alpha(d):=\arccos\left( \frac{ \cos{d}}{1 + \cos{d}}\right)$$ is the angle of the equilateral spherical triangle with side length  $d$.
\item $\sum_{\tau \in I(v)} u_\tau = 2 \pi$ for all vertices $v$ of $G$. Here $I(v)$  is the set of all vertices adjacent edges for a vertex  $v$.
\item If $m_k=3$ then $F_k$ is an equilateral triangle with angles  
 $$u_{k1}=u_{k2}=u_{k3}=\alpha(d).$$
\item  In the case $m_k=4$, $F_k$ is a spherical rhombus and 
 $u_{k1}=u_{k3}$, $u_{k2}=u_{k4}$. Moreover, we have the equality: 
%По правилу Непера для  прямоугольного сферического треугольника получаем равенство: 
$$
\cot{\frac{u_{k1}}{2}}\,\cot{\frac{u_{k2}}{2}} =\cos{d}.
$$

\item In the case $m_k>3$, $F_k=A_1A_2,\ldots,A_{m_k}$ is a convex equilateral spherical polygon with angles $u_{k1},\ldots,u_{m_k}$. The polygon  $F_k$ is uniquely defined (up to isometry) by  its $s:=m_k-3$ angles and  $d$. Then uniquely defined functions $g_i$ and $\zeta_{ij}$ such that $u_{ki}=g_i(u_{k1},\ldots,u_{ks},d)$ and $\dist(A_i,A_j)=\zeta_{i,j}(u_{k1},\ldots,u_{ks},d)$.  from it follows that\\
(a) $u_{ki}=g_i(u_{k1},\ldots,u_{ks},d)$ for $i= m_k-2,m_k-1,m_k$;\\
(b) $\zeta_{i,j}(u_{k1},\ldots,u_{ks},d)\geqslant d$ for $i\ne j$.

\item  Now consider the case when inside  $F_k=A_1A_2,\ldots,A_{m_k}$ there is an isolated vertex.  (It is only if $m_k>5$.) Define
$$
\lambda(u_{k1},\ldots,u_{ks},d):=\max\limits_{p\in F_k}{\min\limits_{i}\{\dist(p,A_i)\}}.
$$
Then $\lambda(u_{k1},\ldots,u_{ks},d)>d.$
\end{enumerate}
\end{prop}

\subsection{Algorithm's description.}
Here we briefly consider our algorithm on enumeration of irreducible contact graphs with $N$ vertices. 
More details can be found in  http://dcs.isa.ru/taras/irreducible/~.

The algorithm consists of two steps:\\
(I) First, we create the list  $L_{N}$ consisting of all graph  with  $N$ vertices that satisfy Corollary 2.1;\\
(II) Using linear approximation of Proposition 4.1 relations from  $L_{N}$ are removed all graphs that cannot be embedded to the sphere. 

\medskip

(I). For the list $L_{N}$ we applied the program {\it plantri}
(see \cite{PLA2})\footnote{Authors of this program are Gunnar Brinkmann and Brendan McKay.}. This program generates non isomorphic planar graphs, including triangulations.  
(In  \cite{PLA1} are given main methods and algorithms of  plantri.)

\medskip

(II). Consider a graph $G$ from $L_{N}$.  We start from the level of approximation  
$\ell=1$. Proposition 4.1 gives possibility to write linear equalities and inequalities for parameters (angles)  $\{u_i\}$ of $G$.

For $\ell=1$ we are using the following relations:\\
(i) $N$ linear equations:  $\sum_{k \in I(v)} u_k = 2 \pi$ (Proposition 4.1 (3).);\\
 (ii) For $\frac{2\pi} {N} \le d $, we obtain  
$ \frac{2\pi}{N} \leqslant \alpha $;\\
(iii)  Proposition 4.1 (5) for a quadrilateral implies  equalities $u_3=u_1, \, u_4=u_2$, and inequalities  $\alpha \leqslant u_i\leqslant 2\alpha , \, i=1,2$;\\
(iv) From the equality $u_2=\rho(u_1,d)$, using   the fact that  $\rho$ is monotonic in both parameters, we obtain maximum and minimum bounds for  $u_2$:
$$ \rho(u_{1,max},d_{min}) \leqslant u_2 \leqslant \rho(u_{1,min},d_{max})$$

\medskip

 So from these linear equalities and inequalities we can obtain maximum and minimum values for each variable.  It gives us a domain $D_1$ that contains all solutions of this system if they there exist. If  $D_1$ is empty, then we can remove $G$ from the list $L_N$. 

This step ``kills'' almost all graphs.

Next we consider $\ell=2$. In this step is divided $D_1$ into two domains and for both we can add the same linear constrains as for $\ell=1$. Moreover, for this step we add new linear constrains for polygons with five and higher vertices.   Some details of this process are given below as well as in our paper \cite[Sec. 4]{MT} and in http://dcs.isa.ru/taras/irreducible/~..

In this level we obtain the parameters domain $D_2$. If this domain is empty, then $G$ cannot be embedded to ${\Bbb S}^2$ and it can be removed from $L_N$.  

 Actually, $\ell=3$ we can repeat previous step, divide  $D_2$ into two domains and obtain additional constrains as for $\ell=2$ for both parts independently. 

We can repeat this procedure more and more times. In fact,  $\ell$ we increase number of subcases. However, practically for every step some subcases are vanished.

We repeat this process for $\ell=1,2,...,m$ and obtain a chain of embedded domains: $$D_m\subset\ldots\subset D_2\subset D_1.$$ If this chain is ended by the empty set, then $G$ can be removed from $L_N$.

In the case if a graph $G$ after certain $m$ steps still ``survived'', i. e.   $D_m\ne\emptyset$,  then it is checked by numerical methods, namely by so called  nonlinear ``solvers''. (We used, in particular, ipopt.) If a solution there exists, then $G$ is declared as a graph that can be embedded, and if not, then $G$ removes from $L_n$.   
%Для этого формируется описание нелинейной задачи на языке AMPL. 
%Если решение найдено успешно, то граф записывается в список найденных
%графов. 
%Для такого графа формируется страница с его описанием (2-$d$
%и 3-$d$ изображения), осуществляектся поиск различных видов решений:
%допустимая точка, решение с минимальным $d$, решение с максимальным
%$d$, внутренняя точка (то есть все углы строго меньше $\pi$).
%Графы, для которых не было найдено решение удаляются из списка проверяемых графов. 
%
%Мы  $N \le 11$ все эти списки оказывались в итоге пустыми.
%То есть для каждого графа было установлено, является он неприводимым
%или нет.

Below we give some details of this algorithm.

\subsection{Linear approximations and the spherical law of cosines.}

In Proposition 4.1 for polygons with four and higher vertices we defined functions $g_i(u_{k1},\ldots,u_{ks},d)$ and   $\zeta_{i,j}(u_{k1},\ldots,u_{ks},d)$. These functions can be calculated by the spherical law of cosines. 
%Эти функции вычисляются на основе сферической теоремы косинусов. Здесь мы рассмотрим более подробно процесс линеаризации этих функций.  

%Мы сделали еще один более глубокий подход к линейному анализу.

%К нашему набору переменных углов всех многоугольников неприводимым графов, мы добавили следующе перменные :

%1. Длина стороны $d$ (которая имеет взаимооднозначное соответствиие углу $a$).

Let $M$ be a polygon with  $m>3$ sides. Let us triangulate $M$ by diagonals and  enumerate the angles of triangles.  
%Разобьем его на треугольники и пронумеруем дополнительные углы.
% и диагонали  как переменные (рис.~\ref{split}). 
For instance, in the case of a pentagon (see Fig. 2) we have nine angles (variables) that  with angles of this pentagon (our variables) are connected by obvious  equations.
%Например, для пятиугольника (Рис. 2) у нас будет 9 углов (переменных), которые с углами самого пятиугольника (наши переменные в алгоритме) связаны очевидными линейными уравнениями. 
%грани $1$ добавлены переменные с именами $F1a1$, \ldots, $F1a9$, $F1d1$, $F1d2$.

%\begin{wrapfigure}{r}{0.2\linewidth}
%\includegraphics[clip,scale=0.5]{pics/irr-5.mps} 
%\label{split}
%\end{wrapfigure}

\begin{figure}[h]
\begin{center}
\includegraphics[clip,scale=0.75]{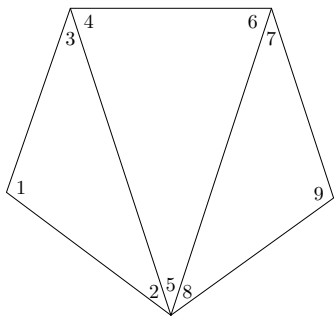}
\end{center}
\caption{Pentagons' angles.}
\label{split}
\end{figure}

%Для разделенных углов грани мы записываем соотношения между старыми и новыми переменными для углов, например ($F1a2+F1a5+F1a8 = V1F1$).

%3. Для каждой переменной мы добавляем дополнительные переменные с префиксами $sin$ и $cos$. Т.е. для переменных  $V5F6$,$d$,$F1a2$ так же есть переменные $cosV5F6$, $sinV5F6$,$cosd$, $sind$, $cosF1a2$, $sinF1a2$.

Actually, if $d$ is fixed, then we have $m-3$ independent variables. To find relations between angles we need just one fact from the spherical trigonometry --- {\em the law of cosines}: 
$$\cos{\phi} = \cos{\theta_1}\cos{\theta_2}+\sin{\theta_1}\sin{\theta_2}\cos\varphi,$$ 
  where for a spherical triangles $ABC$ its sides lengths denoted as $\dist(A,B)=\theta_1, \,
\dist(A,C)=\theta_2, \, \dist(B,C)=\phi$, and $\angle{BAC}=\varphi$. 

 For every triangle from a triangulation we can apply the law of cosines. Then  by interval analysis  all nonlinear inequalities  can be approximated by linear inequalities.     Let us consider some details. 

\subsubsection{Linear inequalities for functions $\sin$ and $\cos$.}

Now we are going to find linear estimations of $f$, where $f(x)=\cos (x)$ or $f(x)=\sin(x)$.
If $x$ lies inside of a given interval  $[x_0 - \delta, x_0 + \delta]$, then 
$$ C \leqslant k x - f(x) \leqslant D.$$

Consider Taylor's series of $f$ at $x_0$:  
$$f(x) = f(x_0) +f'(x_0) (x- x_0) + \ldots + f^{(n)}(x_0)(x- x_0)^n/n! + ...$$

 It is easy to see that the sum of even terms is bounded by  $f''(x_0)(x-x_0)^2/2$,  
and the sum of odd terms started by third  is bounded by  $f'''(x_0)/(x-x_0)^3/6$. 
Therefore, we have 
$k= f'(x_0)$ and
$$C = -f(x_0) + kx_0 + \min(0,f''(x_0)\delta^2/2)-|f'''(x_0)|\delta^3/6,$$
$$D = -f(x_0) + kx_0 + \max(0,f''(x_0)\delta^2/2)+|f'''(x_0)|\delta^3/6,$$

We can substitute by these inequalities all  $\sin$ and $\cos$ functions in equalities and inequalities.

\subsubsection{Linear inequalities for a product $ab$.}
If we have two variables $a$ and $b$, $a \in [a_0 -\delta_a, a_0 + \delta_a]$ and $b \in [b_0 -\delta_b, b_0 + \delta_b]$, then   we have following linear inequalities:
$$ C \leqslant k_aa + k_bb - ab \leqslant D,$$
where $k_a = b_0$, $k_b = a_0$, $C = a_0b_0 -  \delta_a\delta_b$, and $D = a_0b_0 + \delta_a\delta_b$.

\subsubsection{Linear inequalities for $abc$.}
Let 
$a \in [a_0- \delta_a,a_0 + \delta_a]$,
$b \in [b_0- \delta_b,b_0 + \delta_b]$ and 
$c \in [c_0- \delta_c,c_0 + \delta_c]$. Then we have 
$$ C \leqslant k_a a + k_b + k_c c - abc \leqslant D, $$
where 
$k_a = b_0 c_0$, $ k_b = a_0 c_0$, $k_c = a_0 b_0$,
$$C =2a_0b_0c_0 - | a_0 \delta_b \delta_c| - | \delta_a b_0 \delta_c | - | \delta_a \delta_b c_0 | - | \delta_a \delta_b \delta_c|$$
$$D =2a_0b_0c_0 + | a_0 \delta_b \delta_c| + | \delta_a b_0 \delta_c | + | \delta_a \delta_b c_0 | + | \delta_a \delta_b \delta_c|$$

\subsubsection{Linear inequalities for triangles.}
For the law of cosines   
$$
\cos{c}=\cos{a}\cos{b}+\sin{a}\sin{b}\cos{\gamma}.
$$
using intervals for  $\cos$ and $\sin$ we already have linear inequalities. 
For each triangle we write six pairs of linear inequalities for all sides and angles.

\subsection{On optimization of the algorithm.}
It is very important to reduce the time complexity of the algorithm. We are using several ideas for it.

1. For every step $\ell$ using branch-and-bound algorithm we  choose a variable that is divided in two parts.  
For each graph is defined the minimum set of variables that uniquely defined other variables. Then we divide only these variables and $a=\alpha(d)$, that is equal to the angle of a equilateral triangle with side length  $d$ (see Proposition 4.1).
It is essentially  increase speed, in some cases  up to 1000 times.  

2. For variables bounds we have used the following heuristic algorithm. If for some variable we decreased its interval successfully, then we consider  its``neighbors'', i. e. variables that appears together in formulas.

\subsection{On complexity of computations.}
We already noted that there are two steps for enumerating of irreducible contact graphs. In the first step when is created the table  $L_N$ the numbers of graphs graphs grow very fast. For instance, for $N=6,7,8$ , $|L_N|=7,34, 257$. However, for $N=13$: 
$|L_{N}|=94754965$
%Приведем данные о количестве графов  в $L_N$:

In the second step most graphs remove from $L_N$ for  $\ell=1$. However, when  $N$ increases the number of ``bad'' graphs (i. e. graphs that cannot be embedded to the sphere, but survived after many iterations) essentially increases. For these graphs we have to use nonlinear solvers and so computations essentially  increase. That is the main reason why we have tables only up to     $N=11$.

\section{Results}
We applied the method that is discussed above and obtained the following theorem.

\begin{thm} The list of all irreducible contact graphs for  $N=7,8,9,10,11$
on the sphere ${\Bbb S}^2$ is given in tables 5.1--5.5. Here $*$ means that this graph found by Danzer and therefore is D-irreducible, 
and  $**$ means, that this graph is maximal.
 It is also shown in tables bounds for $d$, $d_{min}\leqslant
d \leqslant d_{max}$. (However, note that here values of $d_{min}$
and $d_{max}$ are found numerically and so can be a little different from real.) 
\end{thm}

\subsection{Irreducible graphs with 7 vertices.}
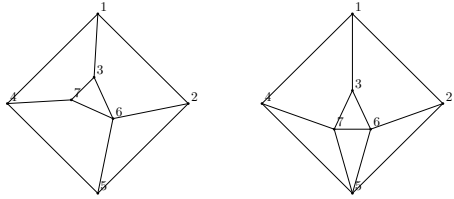
\begin{tabular}{ccc}
$N$ & $d_{min}$ & $d_{max}$ \\
$1*$ & $1.34978$ & $1.35908$ \\
$2**$ & $1.35908$ & $1.35908$ \\
\end{tabular}

\begin{tabular}{cc}
\includegraphics[clip,scale=0.5]{pics/seven1.mps} 
~~~~
\includegraphics[clip,scale=0.5]{pics/seven2.mps}
\end{tabular}

\subsection{Irreducible graphs with 8 vertices.}

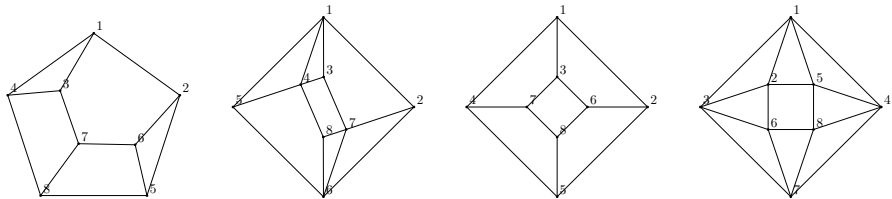
\begin{tabular}{ccc}
$N$ & $d_{min}$ & $d_{max}$ \\
$1$ & $1.17711$ & $1.18349$ \\
$2*$ & $1.28619$ & $1.30653$ \\
$3*$ & $1.23096$ & $1.30653$ \\
$4**$ & $1.30653$ & $1.30653$ \\

\end{tabular}

\begin{tabular}{cc}
\includegraphics[clip,scale=0.5]{pics/eight1.mps} 
~~
\includegraphics[clip,scale=0.5]{pics/eight2.mps}
~~
\includegraphics[clip,scale=0.5]{pics/eight3.mps}
~~
\includegraphics[clip,scale=0.5]{pics/eight4.mps}
\end{tabular}

\subsection{Irreducible graphs with 9 vertices.}
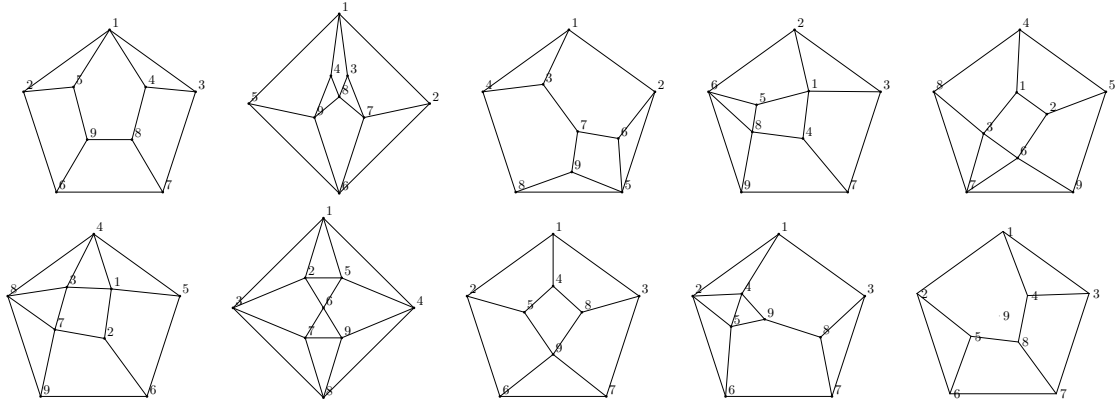
\begin{tabular}{ccc}
$N$ & $d_{min}$ & $d_{max}$ \\
$1$ & $1.14099$ & $1.14143$ \\
$2*$ & $1.22308$ & $1.23096$ \\
$3$ & $1.10525$ & $1.14349$ \\
$4$ & $1.17906$ & $1.18106$ \\
$5$ & $1.15448$ & $1.17906$ \\
$6$ & $1.17906$ & $1.17906$ \\
$7**$ & $1.23096$ & $1.23096$ \\
$8$ & $1.15032$ & $1.18106$ \\
$9*$ & $1.10715$ & $1.14342$ \\
$10$ & $1.17906$ & $1.18428$ \\
\end{tabular}

\begin{tabular}{cc}
\includegraphics[clip,scale=0.5]{pics/nine1.mps} 
~~
\includegraphics[clip,scale=0.5]{pics/nine2.mps} 
~~
\includegraphics[clip,scale=0.5]{pics/nine3.mps} 
~~
\includegraphics[clip,scale=0.5]{pics/nine4.mps} 
~~
\includegraphics[clip,scale=0.5]{pics/nine5.mps} \\
\includegraphics[clip,scale=0.5]{pics/nine6.mps} 
~~
\includegraphics[clip,scale=0.5]{pics/nine7.mps} 
~~
\includegraphics[clip,scale=0.5]{pics/nine8.mps} 
~~
\includegraphics[clip,scale=0.5]{pics/nine9.mps} 
~~
\includegraphics[clip,scale=0.5]{pics/iv_nine1.mps} 
~~
\end{tabular}

\subsection{Irreducible graphs with 10 vertices.}

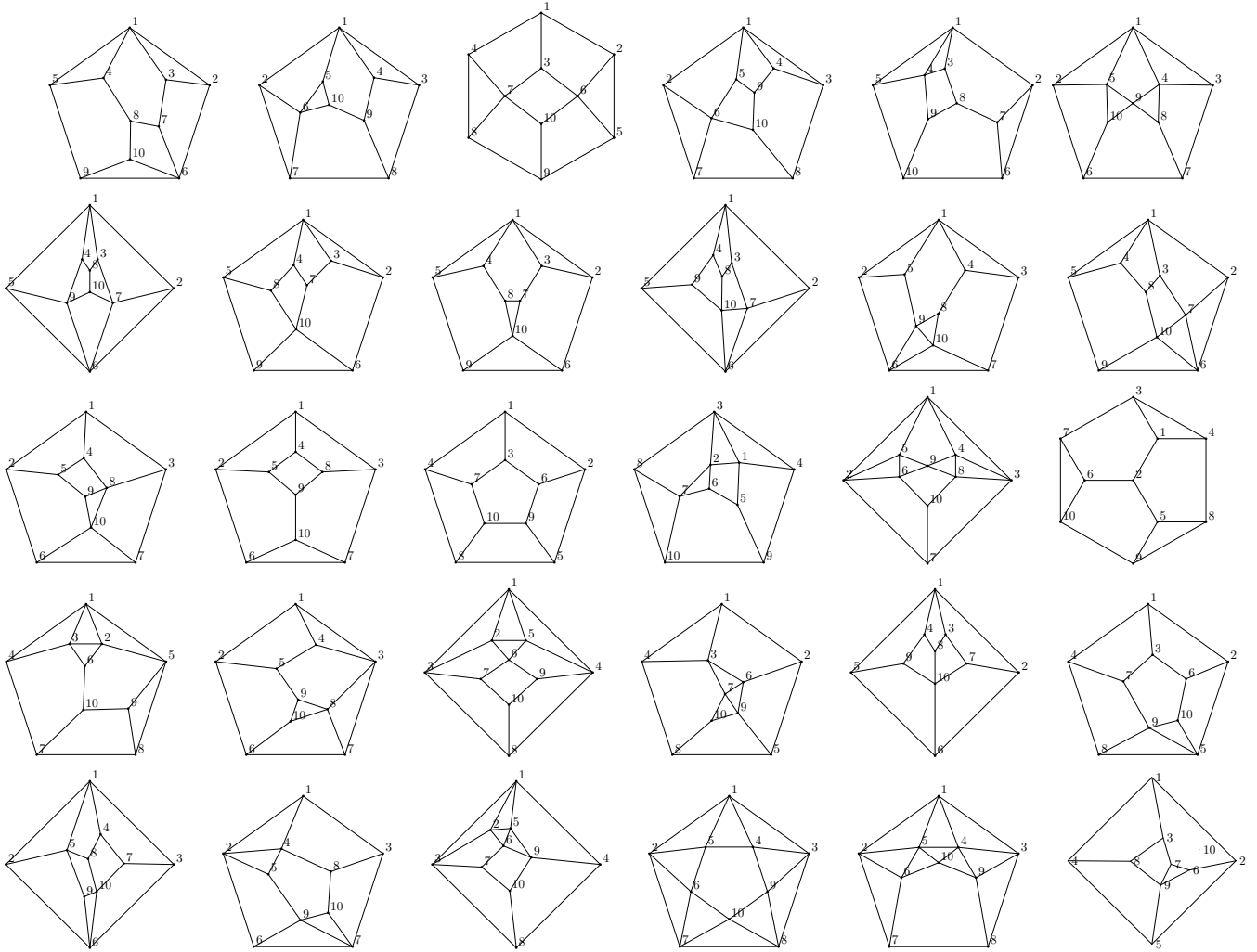
\begin{tabular}{ccccccc}
$N$ & $d_{min}$ & $d_{max}$ & & $N$ & $d_{min}$ & $d_{max}$ \\
$1$ & $1.0839$ & $1.09751$ & & $2$ & $1.08161$ & $1.08439$ \\
$3$ & $1.03067$ & $1.04695$ & & $4$ & $1.10715$ & $1.0988$  \\
$5$ & $1.07529$ & $1.09431$ & & $6$ & $1.09386$ & $1.12285$ \\
$7*$ & $1.15278$ & $1.15448$ & & $8$ & $1.10012$ & $1.10801$ \\
$9$ & $1.06344$ & $1.07834$ & & $10*$ & $1.15074$ & $1.15191$ \\
$11$ & $1.0843$ & $1.08442$ & & $12$ & $1.10055$ & $1.10889$ \\
$13$ & $1.09504$ & $1.10429$ & & $14$ & $1.06032$ & $1.09604$ \\
$15$ & $1.06278$ & $1.1098$ & & $16$ & $1.09567$ & $1.10715$ \\
%\end{tabular}
%\medskip
%\noindent\begin{tabular}{ccccccc}
%$N$ & $d_{min}$ & $d_{max}$ & & $N$ & $d_{min}$ & $d_{max}$ \\
$17**$ & $1.15448$ & $1.15448$ & & $18$ & $0.99865$ & $1.0467$ \\
$19$ & $1.0843$ & $1.0844$ & & $20$ & $1.08334$ & $1.09547$ \\
$21*$ & $1.15341$ & $1.15341$ & & $22$ & $1.0988$ & $1.10608$ \\
$23*$ & $1.14372$ & $1.15191$ & & $24$ & $1.09249$ & $1.1098$ \\
$25*$ & $1.15191$ & $1.15245$ & & $26$ & $1.09658$ & $1.10977$ \\
$27*$ & $1.15191$ & $1.15191$ & & $28*$ & $1.10715$ & $1.10715$ \\
$29*$ & $1.10715$ & $1.10715$ & & $30$ & $1.15103$ & $1.15341$  \\ 
\end{tabular}
\newpage
\begin{tabular}{c}
\includegraphics[clip,scale=0.5]{pics/ten1.mps} 
~~
\includegraphics[clip,scale=0.5]{pics/ten2.mps} 
~~
\includegraphics[clip,scale=0.5]{pics/ten3.mps} 
~~
\includegraphics[clip,scale=0.5]{pics/ten4.mps} 
~~
\includegraphics[clip,scale=0.5]{pics/ten5.mps} 

\includegraphics[clip,scale=0.5]{pics/ten6.mps} 
\end{tabular}
\begin{tabular}{c}
\includegraphics[clip,scale=0.5]{pics/ten7.mps} 
~~
\includegraphics[clip,scale=0.5]{pics/ten8.mps} 
~~
\includegraphics[clip,scale=0.5]{pics/ten9.mps} 
~~
\includegraphics[clip,scale=0.5]{pics/ten10.mps} 
~~
\includegraphics[clip,scale=0.5]{pics/ten11.mps} 
~~
\includegraphics[clip,scale=0.5]{pics/ten12.mps} \\
\end{tabular}
\begin{tabular}{c}
\includegraphics[clip,scale=0.5]{pics/ten13.mps} 
~~
\includegraphics[clip,scale=0.5]{pics/ten14.mps} 
~~
\includegraphics[clip,scale=0.5]{pics/ten15.mps} 
~~
\includegraphics[clip,scale=0.5]{pics/ten16.mps} 
~~
\includegraphics[clip,scale=0.5]{pics/ten17.mps} 
~~
\includegraphics[clip,scale=0.5]{pics/ten18.mps} \\
\end{tabular}
\begin{tabular}{c}
\includegraphics[clip,scale=0.5]{pics/ten19.mps} 
~~
\includegraphics[clip,scale=0.5]{pics/ten20.mps} 
~~
\includegraphics[clip,scale=0.5]{pics/ten21.mps}  
~~
\includegraphics[clip,scale=0.5]{pics/ten22.mps} 
~~
\includegraphics[clip,scale=0.5]{pics/ten23.mps} 
~~
\includegraphics[clip,scale=0.5]{pics/ten24.mps} \\
\end{tabular}
\begin{tabular}{c}
\includegraphics[clip,scale=0.5]{pics/ten25.mps} 
~~
\includegraphics[clip,scale=0.5]{pics/ten26.mps} 
~~
\includegraphics[clip,scale=0.5]{pics/ten27.mps} 
~~
\includegraphics[clip,scale=0.5]{pics/ten28.mps} 
~~
\includegraphics[clip,scale=0.5]{pics/ten29.mps} 
~~
\includegraphics[clip,scale=0.5]{pics/iv_ten1.mps} 
\end{tabular}

\medskip

\subsection{Irreducible graphs with 11 vertices.}
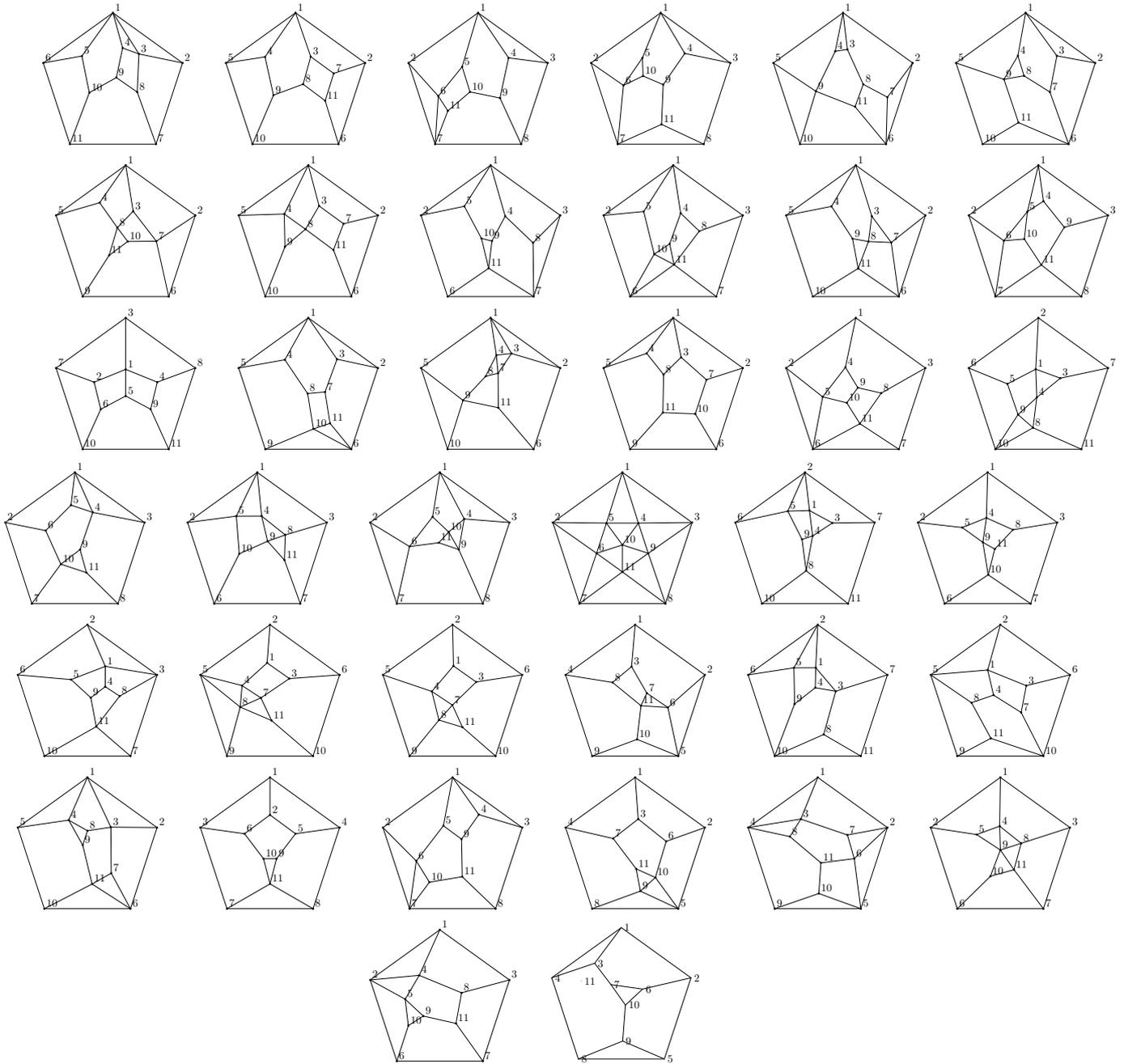
\begin{tabular}{ccccccc}
$N$ & $d_{min}$ & $d_{max}$ & & $N$ & $d_{min}$ & $d_{max}$ \\
$1$ & $1.05601$ & $1.05602$ & & $2$ & $1.0538$ & $1.05842$ \\
$3$ & $1.05834$ & $1.05842$ & & $4$ & $1.04765$ & $1.05455$  \\
$5$ & $1.06975$ & $1.06974$ & & $6$ & $1.06306$ & $1.06308$ \\
$7$ & $1.0522$ & $1.06131$ & & $8$ & $1.06621$ & $1.06846$ \\
$9$ & $1.0538$ & $1.05531$ & & $10$ & $1.0795$ & $1.07961$ \\
$11$ & $1.05331$ & $1.0737$ & & $12$ & $1.07163$ & $1.07197$ \\
$13$ & $1.0404$ & $1.06635$ & & $14$ & $1.04759$ & $1.05637$ \\
$15$ & $1.06974$ & $1.06974$ & & $16$ & $1.02726$ & $1.06117$ \\
$17$ & $1.04712$ & $1.06167$ & & $18$ & $1.06043$ & $1.06209$ \\
\end{tabular}
\begin{tabular}{ccccccc}
$N$ & $d_{min}$ & $d_{max}$ & & $N$ & $d_{min}$ & $d_{max}$ \\
$19$ & $1.05386$ & $1.05947$ & & $20$ & $1.05846$ & $1.05882$ \\
$21$ & $1.0632$ & $1.0636$ & & $22**$ & $1.10715$ & $1.10715$ \\
$23$ & $1.05388$ & $1.06537$ & & $24$ & $1.05375$ & $1.0737$ \\
$25$ & $1.06167$ & $1.0636$ & & $26$ & $1.06506$ & $1.06673$ \\
$27$ & $1.04636$ & $1.05882$ & &  $28$ & $1.05426$ & $1.06822$ \\
$29$ & $1.07832$ & $1.07836$ & & $30$ & $1.07886$ & $1.07962$ \\
$31$ & $1.05429$ & $1.06105$ & & $32$ & $1.00523$ & $1.05671$ \\
$33$ & $1.061$ & $1.06117$ & & $34$ & $1.02751$ & $1.05828$ \\
$35$ & $1.05447$ & $1.06679$ & & $36$ & $1.0561$ & $1.05627$ \\
$37$ & $1.05431$ & $1.05827$ & & $38 (iv)$ & $1.0064$ & $1.03613$ \\
\end{tabular}

\medskip

\begin{tabular}{c}
\includegraphics[clip,scale=0.5]{pics/eleven1.mps} 
~~
\includegraphics[clip,scale=0.5]{pics/eleven2.mps} 
~~
\includegraphics[clip,scale=0.5]{pics/eleven3.mps} 
~~
\includegraphics[clip,scale=0.5]{pics/eleven4.mps} 
~~
\includegraphics[clip,scale=0.5]{pics/eleven5.mps} 
~~
\includegraphics[clip,scale=0.5]{pics/eleven6.mps} \\
~~
\includegraphics[clip,scale=0.5]{pics/eleven7.mps} 
~~
\includegraphics[clip,scale=0.5]{pics/eleven8.mps} 
~~
\includegraphics[clip,scale=0.5]{pics/eleven9.mps} 
~~
\includegraphics[clip,scale=0.5]{pics/eleven10.mps} 
~~
\includegraphics[clip,scale=0.5]{pics/eleven11.mps} 
~~
\includegraphics[clip,scale=0.5]{pics/eleven12.mps} \\
~~
\includegraphics[clip,scale=0.5]{pics/eleven13.mps} 
~~
\includegraphics[clip,scale=0.5]{pics/eleven14.mps} 
~~
\includegraphics[clip,scale=0.5]{pics/eleven15.mps} 
~~
\includegraphics[clip,scale=0.5]{pics/eleven16.mps} 
~~
\includegraphics[clip,scale=0.5]{pics/eleven17.mps} 
~~
\includegraphics[clip,scale=0.5]{pics/eleven18.mps} \\
\end{tabular}
\begin{tabular}{c}
\includegraphics[clip,scale=0.5]{pics/eleven19.mps} 
~~
\includegraphics[clip,scale=0.5]{pics/eleven20.mps} 
~~
\includegraphics[clip,scale=0.5]{pics/eleven21.mps} 
~~
\includegraphics[clip,scale=0.5]{pics/eleven22.mps} 
~~
\includegraphics[clip,scale=0.5]{pics/eleven23.mps} 
~~
\includegraphics[clip,scale=0.5]{pics/eleven24.mps} \\
~~
\includegraphics[clip,scale=0.5]{pics/eleven25.mps} 
~~
\includegraphics[clip,scale=0.5]{pics/eleven26.mps} 
~~
\includegraphics[clip,scale=0.5]{pics/eleven27.mps} 
~~
\includegraphics[clip,scale=0.5]{pics/eleven28.mps} 
~~
\includegraphics[clip,scale=0.5]{pics/eleven29.mps} 
~~
\includegraphics[clip,scale=0.5]{pics/eleven30.mps} \\
~~
\includegraphics[clip,scale=0.5]{pics/eleven31.mps} 
~~
\includegraphics[clip,scale=0.5]{pics/eleven32.mps} 
~~
\includegraphics[clip,scale=0.5]{pics/eleven33.mps} 
~~
\includegraphics[clip,scale=0.5]{pics/eleven34.mps} 
~~
\includegraphics[clip,scale=0.5]{pics/eleven35.mps} 
~~
\includegraphics[clip,scale=0.5]{pics/eleven36.mps} \\
~~
\includegraphics[clip,scale=0.5]{pics/eleven37.mps} 
~~
\includegraphics[clip,scale=0.5]{pics/iv_eleven1.mps} 
~~
\end{tabular}

\medskip

\medskip

\medskip

\medskip

\medskip

\medskip

O. R. Musin, IITP RAS и UTB  (University of Texas at
Brownsville).

 {\it E-mail:} oleg.musin@utb.edu

\medskip

A. S. Tarasov, IITP RAS

{\it E-mail:} tarasov.alexey@gmail.com


\begin{thebibliography}{99}

\bibitem{AZ}
M. Aigner and G.M. Ziegler, Proofs from THE BOOK, Springer, 1998
(first ed.) and 2002 (second ed.)
\bibitem{Ans}
K. Anstreicher, The thirteen spheres: A new proof, Discrete Comput.
Geom. {\bf 31}(2004), 613-625.

%\bibitem{BV}
%C. Bachoc and F. Vallentin, New upper bounds for kissing numbers
%from semidefinite programming,
%J. Amer. Math. Soc. {\bf 21} (2008), 909-924.

\bibitem{BM}
A. Barg and O.~R. Musin, Codes in spherical caps,
Advances in Mathematics of Communication, {\bf 1} (2007), 131-149.

\bibitem{Bor11}
K. B\"or\"oczky, The problem of Tammes for $n = 11$, Studia.
Sci. Math. Hungar. {\bf 18} (1983) 165-171.

\bibitem{Bor}
K. B\"or\"oczky, The Newton-Gregory problem revisited, In: Discrete
Geometry, A. Bezdek (ed.), Dekker, 2003,  103-110.

\bibitem{BS}
K. B\"or\"oczky, L. Szab\'o, Arrangements of 13 points on a sphere,
In: Discrete Geometry, A. Bezdek (ed.), Dekker, 2003, 111-184.

\bibitem{BS14}
K. B\"or\"oczky, L. Szab\'o, Arrangements of 14, 15, 16 and 17
points on a sphere, Studi. Sci. Math. Hung. {\bf 40} (2003),
407-421.

\bibitem{BDM}
P. Boyvalenkov, S. Dodunekov and O. R. Musin, A survey on the kissing numbers, Serdica Mathematical Journal,  {\bf 38} (2012), 507-522.

\bibitem{BMP}
P. Brass, W.O.J. Moser, J. Pach, Research problems in discrete geometry, Springer-Verlag, 2005.

\bibitem{PLA1}
G. Brinkmann and B. D. McKay, Fast generation of planar graphs (expanded edition), http://cs.anu.edu.au/~bdm/papers/plantri-full.pdf

%\bibitem{Cas}
%B. Casselman, The difficulties of kissing in three dimensions,
%Notices Amer. Math. Soc., {\bf 51}(2004), 884-885.

\bibitem{CS}
J.H. Conway and N.J.A. Sloane, Sphere Packings, Lattices, and Groups, New York, Springer-Verlag, 1999 (Third Edition).

%\bibitem{ConRig}
%R. Connelly, Generic global rigidity, Discrete Comput. Geom. {\bf 33} (2005), no.4, 549-563.

\bibitem{Appl04}
A. Donev, S. Torquato, F. H. Stillinger, and R. Connelly, Jamming in Hard Sphere and Disk Packings, Journal of Applied Physics, 95, 989-999 (2004). 

\bibitem{Dan}
L. Danzer, Finite point-sets on ${\bf S}^2$ with minimum distance
as large as possible, Discr. Math., {\bf 60} (1986), 3-66.


\bibitem{FeT0}
L. Fejes T\'oth, \"Uber die Absch\"atzung des k\"urzesten Abstandes
zweier Punkte eines auf einer Kugelfl\"ache liegenden Punktsystems,
Jber. Deutch. Math. Verein. {\bf 53} (1943), 66-68.


\bibitem{FeT}
L. Fejes T\'oth, Lagerungen in der Ebene, auf der Kugel und in
Raum, Springer-Verlag, 1953; Russian translation, Moscow, 1958


%\bibitem{Hales}
%T. Hales, The status of the Kepler conjecture, Mathematical Intelligencer
%{\bf 16}(1994), 47-58.

%\bibitem{Hop}
%R. Hoppe, Bemerkung der Redaction, Archiv Math. Physik (Grunet)
%{\bf 56} (1874), 307-312.

\bibitem{HabvdW}
W. Habicht und B.L. van der Waerden, Lagerungen von Punkten auf der Kugel, Math. Ann. {\bf 123} (1951),  223-234.

\bibitem {Appl10} A. B. Hopkins, F. H. Stillinger, and S. Torquato, Densest Local Sphere-Packing Diversity: General concepts and application to two dimensions, Physical Review E, 81, 041305 (2010). 

\bibitem{Hs}
W.-Y. Hsiang, Least action principle of crystal formation of
dense packing type and Kepler's conjecture, World Scientific,
2001.

\bibitem{Lee}
J. Leech, The problem of the thirteen spheres, Math. Gazette
{\bf 41} (1956), 22-23.


%V.I. Levenshtein, On bounds for packing in $n$-dimensional Euclidean
%space, Sov. Math. Dokl.
%{\bf 20}(2), 1979, 417-421.

\bibitem{Ma}
H. Maehara, Isoperimetric theorem for spherical polygons and
the problem of 13 spheres, Ryukyu Math. J.,
{\bf 14} (2001), 41-57.

\bibitem{Manew}
H. Maehara, The problem of thirteen spheres - a proof for undergraduates,
European Journal of Combinatorics, {\bf 28} (2007), 1770-1778.

\bibitem{Mus1}
O.R. Musin. The problem of the twenty-five spheres// Russian Math. Surveys,  {\bf 58} (2003), 794-795. 

\bibitem{Mus13}
O.~R. Musin, The kissing problem in three dimensions, Discrete
Comput. Geom., {\bf 35} (2006), 375-384.


\bibitem{Mus3}
O.~R. Musin, The one-sided kissing number in four dimensions,
Periodica Math. Hungar., {\bf 53} (2006), 209-225.

\bibitem{Mus4}
O.~R. Musin, The kissing number in four dimensions, Ann. of Math., {\bf 168} (2008),  1-32.

\bibitem{Mus5}
O.~R. Musin, Bounds for codes by semidefinite programming, Proc. Steklov Inst. Math. {\bf 263} (2008), 134-149.

\bibitem{Mus6}
O.R. Musin, Positive definite functions in distance geometry, European Congress of Mathematics Amsterdam, 14-18 July, 2008, 115-134, EMS Publ. 2010.

\bibitem{MT}
O.~R. Musin and A.~S. Tarasov, The Strong Thirteen Spheres Problem, Discrete \& Comput. Geom., {\bf 48} (2012), 128-141.


%\bibitem{OdS}
%A.M. Odlyzko and N.J.A. Sloane, New bounds on the number of unit
%spheres that
%can touch a unit sphere in $n$ dimensions, J. of Combinatorial
%Theory A, {\bf 26} (1979), 210-214.

\bibitem{PLA2}
plantri and fullgen, http://cs.anu.edu.au/~bdm/plantri/

%\bibitem{PZ}
%F. Pfender and G.M. Ziegler, Kissing numbers, sphere packings,
%and some unexpected proofs, Notices Amer. Math. Soc., {\bf 51}(2004),
%873-883.

\bibitem{Rob}
R.M. Robinson, Arrangement of 24 circles on a sphere, Math. Ann.
{\bf 144} (1961), 17-48.


\bibitem{SvdW1}
K. Sch\"utte and B.L. v. d. Waerden, Auf welcher Kugel haben
5,6,7,8 oder 9
Punkte mit Mindestabstand 1 Platz? Math. Ann. {\bf 123} (1951),
96-124.

\bibitem{SvdW2}
K. Sch\"utte and B.L. van der Waerden, Das Problem der dreizehn
Kugeln, Math. Ann. {\bf 125} (1953), 325-334.

%\bibitem{Sz1}
%G.G. Szpiro, Kepler's conjecture, Wiley, 2002.



\bibitem{Tam}
 R.M.L. Tammes, On the Origin Number and Arrangement of the Places
of Exits on the Surface of Pollengrains, Rec. Trv. Bot. Neerl.
{\bf 27} (1930), 1-84.


\bibitem{vdW}
B.L. van der Waerden, Punkte auf der Kugel. Drei Zus\"atze, Math. Ann. {\bf 125} (1952) 213-222.

%R. Connely, Discrete Comput. Geom., Volume 33, Number 4, April 2005, pages 549-563


\end{thebibliography}
\end{document}